\documentclass[10pt,a4paper,oneside,leqno]{article}
\usepackage{amsfonts,amssymb}
\usepackage{amscd}
\usepackage{color}

\newtheorem{defn}{Definition}

\newtheorem{lemma}{Lemma}
\newtheorem{theoremn}{Theorem}
\newtheorem{observen}{Observation}

\newcommand{\fnomial}[2]{ {{#1} \choose {#2}}_F }

\vspace{0.5cm}

 4
\font\ebf=cmbx8
\font\erm=cmr8

\usepackage{graphicx, graphics}

\begin{document}
\begin{center}
\noindent {\Large \textsc{On cobweb posets' most relevant codings}} \textcolor{blue}{(\textbf{***})}  \\ 
\vspace{0.5cm}
\noindent A. Krzysztof Kwa\'sniewski (*),  M. Dziemia\'nczuk (**) \\
\vspace{0.2cm}
\noindent {\erm (*)  member of the Institute of Combinatorics and its Applications  }\\
\noindent {\erm PL-15-021 Bia\l ystok, ul. Kamienna 17, Poland}\\
\noindent {\erm  High School of Mathematics and Applied Informatics}\\
\noindent {\erm e-mail: kwandr@gmail.com}\\
\noindent {\erm (**) Institute of Informatics, University of Gda\'nsk}\\
\noindent {\erm PL-80-952 Gda\'nsk, st. Wita Stwosza 57, Poland}\\
\noindent {\erm e-mail: mdziemianczuk@gmail.com} 
\end{center}
\vspace{0.4cm}

\noindent {\ebf SUMMARY:}

\vspace{0.1cm}

\noindent {\small  One considers here acyclic digraphs named   KoDAGs \textbf{(****)} which represent  the outmost general chains  of di-bi-cliques denoting thus the outmost general chains of  binary relations. Because  of this fact KoDAGs  start to become an outstanding concept of nowadays investigation. We propose here examples of  codings of KoDAGs looked upon as infinite hyper-boxes as well as  chains of rectangular  hyper-boxes in $N^\infty$. Neither of KoDAGs' codings considered here is a  poset isomorphism  with $\Pi = \langle P, \leq\rangle$. Nevertheless every example of  coding supplies a new view on possible investigation of KoDAGs properties.  The codes proposed here down are by now recognized as most relevant codings  for practical purposes including visualization. 
\noindent More than that.  Employing quite arbitrary sequences $F=\{n_F\}_{n\geq 0}$ infinitely many  new representations of  natural numbers called an  $F$- base or  base-$F$ number system representations are introduced. These constitute mixed radix-type  numeral systems. $F$- base non-standard positional numeral systems in which the numerical base varies from position to position have picturesque interpretation due to KoDAGs graphs and their correspondent posets which in turn  are  endowed on their own with combinatorial interpretation of uniquely assigned to  KoDAGs $F-nomial$ coefficients.  The base-$F$ number systems are used for KoDAGs' coding and are interpreted as chain coordinatization in KoDAGs pictures as well as systems of  infinite number of  boxes' sequences of $F$-varying containers capacity of subsequent boxes.  Needless to say how crucial is this base-F number system for KoDAGs - hence - consequently for arbitrary chains of binary relations. New $F$-based numeral systems are umbral base-$F$ number systems in a sense to be explained in what follows.} 

\vspace{0.2cm}

\noindent  \textbf{(****)} - for the history of this name see The Internet Gian-Carlo Polish Seminar \textbf{Subject 1} \textit{oDAGs and KoDAGs in Company} - here $\Rightarrow$ \emph{http://ii.uwb.edu.pl/akk/sem/sem\_rota.htm}.

\vspace{0.2cm}

\noindent  \textcolor{blue}{(\textbf{***})} This is The Internet Gian-Carlo Polish Seminar article, update of [0]\\
No \textbf{2}, \textbf{Subject 2}, \textbf{2009-02-25}\\
\noindent \emph{http://ii.uwb.edu.pl/akk/sem/sem\_rota.htm}\\
 
\vspace{0.2cm}


\vspace{0.1cm}

\noindent AMS Classification Numbers: 06A07 ,05C78, 11A63,
\vspace{0.1cm}

\vspace{0.2cm}
\noindent Affiliated to The Internet Gian-Carlo Polish Seminar: \\
\noindent \emph{http://ii.uwb.edu.pl/akk/sem/sem\_rota.htm}

\vspace{0.4cm}

\section{Introduction. Cobweb poset definition.}
\noindent Let any natural numbers valued sequence $F$-sequence be chosen.  (Sometimes  value zero might be admitted as for example $F_0 = 0$ in the case of $F$ being Fibonacci sequence). Then $F$ designates all structures considered here. Here  a directed acyclic graph, also called DAG, is a directed graph with no directed cycles.   Note that  posets  $\Pi = \langle P, \leq\rangle$ defined below are not lattices except for trivial case.

\subsection{Plane $N_0\times N$ definition of $\Pi = \langle P, \leq\rangle$ \cite{12,19,18,17,27}}

Cobweb poset $\Pi = \langle P, \leq\rangle$ is defined via its Hasse digraph as follows:

$$
	P = \bigcup_{s\geq 0}{\Phi_s}, \ \ \ \ \Phi = \left\{ \langle j, s \rangle : 1 \leq j \leq s_F \right\}
$$

\noindent $s$-level, $j$ indicates position of a point in $s$-th level

$$
	 \langle j, s \rangle \leq  \langle k, t \rangle  \Leftrightarrow  s < t \vee (j = k \wedge s = t).
$$

\noindent Any cobweb poset  and  any cobweb subposet  graph  $P_k = \bigcup_{s\geq 0}^{k}{\Phi_s}$ is a \textbf{DAG} \cite{19,18,17}.

\vspace{0.2cm}
\noindent More than that - all these posets as well  as  any  cobweb layer  $\langle \Phi_k \rightarrow \Phi_n \rangle = \bigcup_{s\geq k}^{n}{\Phi_s} $ are faithfully represented via  orderable DAGs  i.e. \textbf{oDAGs} \cite{28,6,29} and since recently these have being called KoDAGs  \cite{1,2,6,7,8,9}.

\vspace{0.2cm}
\noindent While cobweb layer  $\langle \Phi_k \rightarrow \Phi_n \rangle$, $n\geq k$ is an  $(n-k+1)$-level digraph oDAG and cobweb subposet $P_n$ is an  $(n+1)$-level digraph specifically $\langle \Phi_k \rightarrow \Phi_{k+1} \rangle$ is just a complete bipartite digraph

\begin{equation}\label{eq:1}
	\langle \Phi_k \rightarrow \Phi_{k+1} \rangle = K_{F_k, F_{k+1}}, \ \ \ \ k = 0, 1, 2, ...
\end{equation}

\subsection{Chain of di-bi-cliques' definition of $\Pi = \langle P, \leq \rangle$  \cite{19}}

Identification (\ref{eq:1}) might be considered as source of the definition  of  $F_0$  rooted $F$-cobweb graph $P(F) = P$. The usual convention is to  establish  $F_0 = 1$.

\vspace{0.4cm}
\noindent \textbf{Namely:} 
Cobweb poset $\Pi = \langle P, \leq \rangle$ is defined via its  Hasse digraph as a \textbf{chain of complete binary relations} each one link of the chain being represented by by its di-bi-clique  graph  i.e. complete \textbf{bipartite digraph} $K_{F_k, F_{k+1}}$; $k=0,1,2,...$ where $\langle j,s \rangle \leq \langle k,t \rangle \Leftrightarrow  s < t \vee j = k \wedge s = t$.

\vspace{0.2cm}
\noindent \textbf{Obviously} we thus establish natural and obvious bijection uniquely coding \emph{complete} relations' chain via its KoDAG.

\vspace{0.2cm}
\noindent \textbf{Obviously}  any chain of binary relations is obtainable from the cobweb poset chain of  complete relations via deleting arcs in di-bicliques of the complete relations chain as of course any relation  $R_k$ as a subset of  $\Phi_k \times \Phi_{k+1}$ is represented by a  bipartite digraph  $D_k$. Recall then that a "complete relation"  $C_k$ is identified by definition with its  di-biclique graph  $K_{F_k, F_{k+1}}$. A digraph  $D_k$ is a sub- digraph of  $C_k$.

\vspace{0.2cm}
\noindent \textbf{In brief then} \cite{19} \emph{cobweb posets' and thus KoDAGs's defining  di-bicliques are  links  of any   complete relations' chain.} 

\vspace{0.2cm}
\noindent For more on that  see \cite{19}  from where we quote

\begin{quotation}
DAGs  considered  as a generalization of trees  have a lot of   applications in computer science, bioinformatics,  physics and many natural activities of humanity and nature.  For example in information categorization systems, such as folders in a computer or in Serializability Theory of Transaction Processing Systems and many others. Here we introduce specific DAGs as generalization of trees being inspired by algorithm of the Fibonacci tree growth. For any given natural numbers valued sequence the graded (layered) cobweb posets` DAGs  are equivalently representations of a chain of binary relations. Every relation of the cobweb poset chain is biunivocally represented by the uniquely designated  \textbf{complete} bipartite digraph-a digraph which is a di-biclique  designated  by the very  given sequence. The cobweb poset is then to be identified with a chain of di-bicliques i.e. by definition - a chain of complete bipartite one direction digraphs.   Any chain of relations is therefore obtainable from the cobweb poset chainof complete relations  via deleting  arcs (arrows) in di-bicliques.\\                                                                                                           
Let us underline it again : \textit{any chain of relations is obtainable from the cobweb poset chain of  complete relations via deleting arcs in di-bicliques of the complete relations chain.} For that to see note that any relation  $R_k$ as a subset of  $A_k \times A_{k+1}$ is represented by a  one-direction bipartite digraph  $D_k$.  A "complete relation"  $C_k$ by definition is identified with its one direction di-biclique graph $d-B_k$.  Any  $R_k$ is a subset of  $C_k$. Correspondingly one direction digraph  $D_k$ is a subgraph of an one direction digraph of $d-B_k$.\\                       The one direction digraph of  $d-B_k$ is called since now on  \textbf{the di-biclique }i.e. by definition - a complete bipartite one direction digraph.   Another words: cobweb poset defining di-bicliques are links of a complete relations' chain. (end of quote)
\end{quotation}

\begin{figure}[ht]
\begin{center}
	\includegraphics[width=100mm]{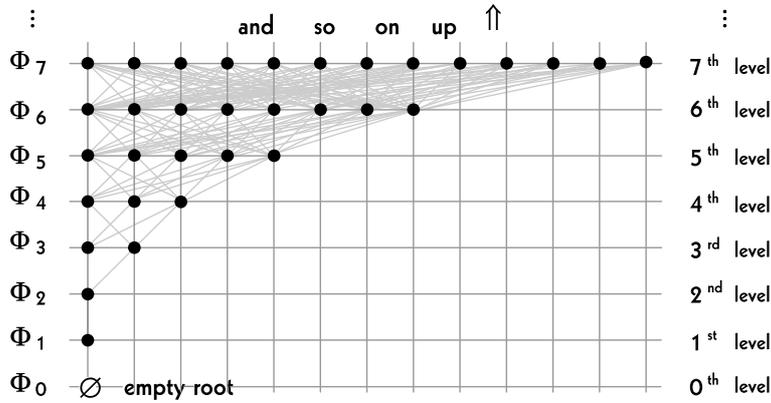}
	\caption{The $s$-th level in $\mathbb{N}\times \mathbb{N}\cup\{0\}$ \label{fig:1}}
\end{center}
\end{figure}

\begin{figure}[ht]
\begin{center}
	\includegraphics[width=75mm]{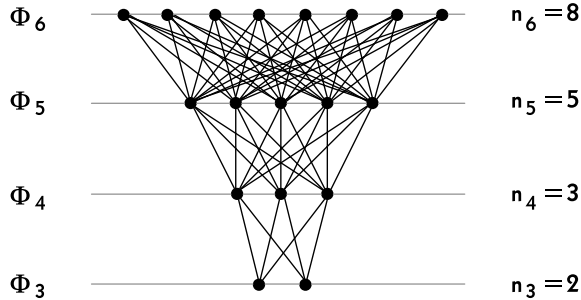}
	\caption{Display of four levels of Fibonacci numbers' finite Cobweb sub-poset \label{fig:2}}
\end{center}
\end{figure}

\begin{figure}[ht]
\begin{center}
	\includegraphics[width=75mm]{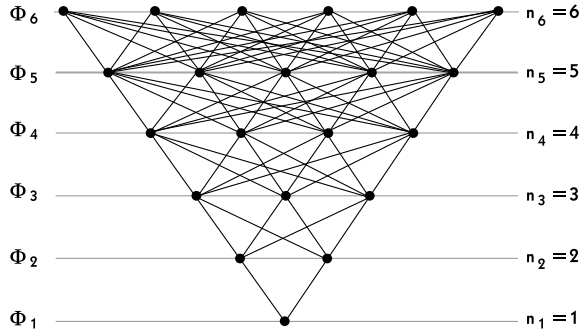}
	\caption{Display of Natural numbers' finite prime Cobweb poset \label{fig:3}}
\end{center}
\end{figure}

\begin{figure}[ht]
\begin{center}
	\includegraphics[width=100mm]{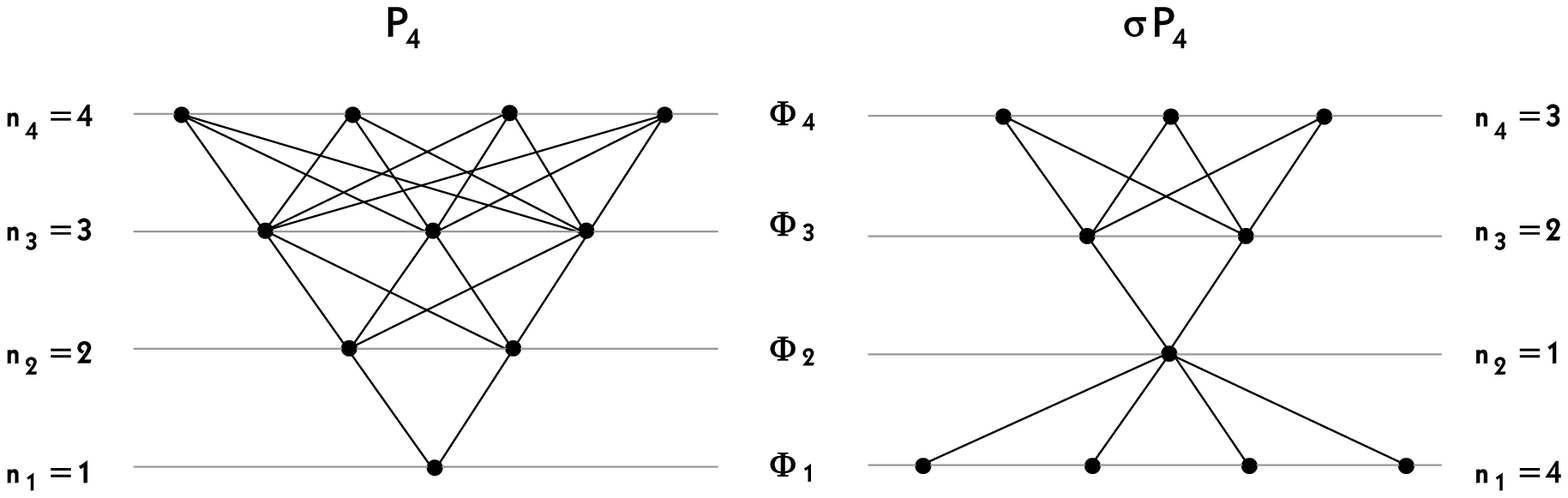}
	\caption{Display of block $\sigma P_m$ obtained from $P_m$ and permutation $\sigma$ \label{fig:4}}
\end{center}
\end{figure}

\noindent Here come some examples of KoDAGs (Fig. \ref{fig:1} - \ref{fig:4}).

\noindent Notation used in above figure an figures to follow comes from \cite{19,18,17} and \cite{27} where

$$
	\sigma P_m = C_m[F; \sigma <F_1, F_2,...,F_m>]
$$  is an equipotent sub-poset obtained from $P_m$ with help of a permutation $\sigma$  of the sequence $<F_1, F_2,...,F_m>$ encoding  $m$ layers  of $P_m$ thus obtaining the equinumerous sub-poset $\sigma P_m$   with the sequence $\sigma <F_1, F_2,...,F_m>$ encoding  now $m$ layers of $\sigma P_m$. Then $P_m = C_m[F; <F_1, F_2,...,F_m>].$

\vspace{0.2cm}
\noindent Further readings: \cite{12,19,18,17,25,1,2,3,6,7,8,9,27}.

\vspace{0.2cm}
\noindent \textbf{Motto for further investigations:}

\vspace{0.2cm}
\emph{ KoDAGs represent the outmost general chains of di-bi-cliques denoting thus the outmost general chains     of  binary relations. Because  that fact KoDAGs start to become  an outstanding concept of nowadays investigations.}

\section{Cobweb posets' coding via $N^\infty$ lattice boxes and more}

\textbf{Basic:} $N^\infty = N\times N\times N \times ... = \times_{s\geq 0}{N_s}$, $N_s \equiv N$. Neither of codings introduced in this section  is the  posets' isomorphism  with $\langle \Pi, \leq\rangle$.

\vspace{0.2cm}
\noindent Recall that posets $\Pi = \langle P, \leq \rangle$ defined below are not lattices except for trivial  case.

\noindent At the same time for the sake of Kwa\'sniewski combinatorial interpretation of cobweb poset's inherited $F$-nomial coefficients \cite{19,10,11,22,23,12,18,17,4,15,14,13,5} - \textbf{maximal chains}  of layers $\langle \Phi_k \rightarrow \Phi_n \rangle$ were and are looked upon by  Kwa\'sniewski  expresis verbis as the \textbf{set points} of the  set

$$
	\big\{ \mathrm{maximal\ chains\ in\ } \langle \Phi_k \rightarrow \Phi_n \rangle \big\} \equiv
	C_{max}\big( \langle \Phi_k \rightarrow \Phi_n \rangle \big) \equiv
	C^{k,n}_{max}.
$$

\vspace{0.2cm}
\noindent On the other hand - Dziemia\'nczuk [3,2,1] - while investing tiling problem  posed in \cite{19,18,17}
had observed  in practice of computer calculations an usefulness of his so called  geometric coding  of the posets $\Pi = \langle P, \leq \rangle$. In the course of further joint searching his concept clarified to  be of much more than of  technical  relevance.  Let us present it as it sounds now taking into account  both combinatorial and  \emph{geometric-tiling}  points of viewing posets  $\Pi = \langle P, \leq \rangle$   simultaneously.

\subsection{Geometric $N^\infty$ coding of $\Pi$ via the set $V$ of maximal chains of $\Pi$ }

The poset $K$  is to be defined  as  $K \equiv \langle V, \leq_{\Pi} \rangle$, where $V \subset N^\infty$ is selected to be the set of points from $N^\infty$ forming the infinite  discrete  $F$-hyper-box   (discrete rectangular hyper-box designated by $F$-sequence...  or in everyday parlance just hyper-box  i.e.
$V = [0_F]\times [1_F]\times [2_F] \times ... $ where $[0_F]$ might correspond as in Fibonacci sequence case  to an "empty root" and just by convention has one point $\{\emptyset \}$.
$$
{
	V \equiv \big\{ \chi = \left(c_0,c_1,c_2... \right) : c_n\in \{0,1,...,n_F-1\}, n\!\geq\!0 \big\} \equiv
\atop
	\equiv \big\{ \chi = \left( \phi_0,\phi_1,\phi_2... \right) : \phi_n \in \Phi_n, n\!\geq\!0 \big\}
}
$$

\noindent or in plane $N_0 \times N$ presentation 

$$
	V \Leftrightarrow \big\{ \chi = \left( \langle c_0,0\rangle, \langle c_1,1\rangle,... \right) :1\leq c_n \leq n_F, n \geq 0 \big\}
$$

\vspace{0.2cm}
\noindent with the point $\chi = \left(c_0,c_1,c_2,...\right)$, $\chi \in V$ interpreted as a \emph{maximal-chain} in $\Pi$. $V$ denotes  the family  of maximal chains of  $\Pi$.

\vspace{0.2cm}
\noindent \textbf{Example:} for $\leq_\Pi = \leq_V$ where for $\chi_1 = (c_0,c_1,c_2,...)$, $\chi_2 = (d_0,d_1,d_2,...)$ we define

$$
	\chi_1 \leq_V \chi_2   \Leftrightarrow  \big( \exists_{n\in N}\ c_n < d_n \wedge \forall_{m>n}\ c_m \leq d_m \big)
$$

\vspace{0.2cm}
\noindent le voila a chain of coded $\Pi$ poset maximal chains \\
$(0,0,0,...) \leq_V (0,0,1,0,0,...) \leq_V (0,0,0,1,0,...) \leq_V (0,0,1,1,0,...)\ ... \leq_V ...$

\vspace{0.2cm}
\noindent The above is an example of a chain in $K$. It is characteristic example as $\leq_V$ is a total order.

\vspace{0.2cm}
\noindent Consider now $V^* = \left\{ (c_0,c_1,...,c_m,0,0,...)\in V \right\}$. It appears to be a chain of chains.

\begin{lemma}\label{lem:1}
$\leq_V$ is total order in $V^*$ hence $\langle V^*, \leq_V \rangle$ is a chain.
\end{lemma}

\noindent \textit{Proof.} Indeed. Let $\chi_1 \leq_V \chi_2$, $\chi_1 \neq \chi_2$. \\
Then either 
$$
	\chi_1 \leq_V \chi_2 \Leftrightarrow \left(\exists n\in N : c_n < d_n \wedge \forall m>n : c_m \leq d_m\right)
$$
 or (note that $\chi_1 \neq \chi_2$ is true) \\
$\neg (\chi_1 \leq_V \chi_2)$ 
$\Leftrightarrow$ $\neg \left(\exists n\in N : c_n < d_n \wedge \exists n\in N, \forall m>n : c_m \leq d_m\right)$ \\
$\Leftrightarrow$ $\neg (\exists n\in N : c_n < d_n)$ $\vee$ $\neg (\exists n \in N, \forall m > n : c_m \leq d_m )$ \\
$\Leftrightarrow$ $(\forall n\in N : d_n \leq c_n)$ $\vee$ $\neg (\exists n \in N, \forall m > n : c_m \leq d_m )$ \\
$\Leftrightarrow$ $\big[ (\forall n\in N : d_n \leq c_n) \wedge \chi_1 \neq \chi_2 \big]$ 
	$\vee$ $\neg (\exists n \in N, \forall m > n : c_m \leq d_m )$ $\vee$ $(\exists n,m\in N, m>n : d_m < c_m)$\\
$\Leftrightarrow$ $\chi_2 \leq_V \chi_1$ 
	$\vee$ $\big[ \neg (\exists n \in N, \forall m > n : c_m \leq d_m ) \vee (\exists n,m\in N, m>n : d_m < c_m) \big]$ \\
$\Leftrightarrow$ $\chi_2 \leq_V \chi_1$ 
	$\vee$ $\big[ \forall m\in N : d_m < c_m \vee (\exists n,m\in N: m>n : d_m < c_m) \big]$ \\
$\Leftrightarrow$  $\chi_2 \leq_V \chi_1$ $\vee$  $\chi_2 \leq_V \chi_1$
	$\vee$ $(\exists n,m\in N : m>n : d_m < c_m)$ \\
$\Leftrightarrow$  $\chi_2 \leq_V \chi_1$ $\vee$  $\chi_2 \leq_V \chi_1$ $\vee$  $\chi_2 \leq_V \chi_1$ \\
$\Leftrightarrow$  $\chi_2 \leq_V \chi_1$ \\
Let $\chi_1 \neq \chi_2$. Then either $\chi_1 \leq_V \chi_2$ or $\chi_2 \leq_V \chi_1$ $\blacksquare$

\subsection{Natural Partial Order Relations $\leq_\Pi$ defining coding poset $K \equiv \langle V, \leq_\Pi \rangle$ }

\noindent \textbf{2.2.1} \  $\leq_V$, $\Gamma \equiv \langle V, \leq_V \rangle$

\vspace{0.2cm}
\noindent For the reasons to be apparent soon (see also Section 4 pictures) the codings to follow shall  be called discrete \textbf{geometric}.

\vspace{0.2cm}
\noindent A partial  order relation $\leq_\Pi = \leq_V$  already introduced and used for the sake of computer calculations 
with maximal chains' counting  by Dziemia\'nczuk 

\noindent [\emph{http://www.dejaview.cad.pl/cobwebposets.html}]  has been defined as follows

$$
	\chi_1 \leq_V \chi_2   \Leftrightarrow  \big( \exists_{n\in N}\ c_n < d_n \wedge \forall_{N \ni m>n}\ c_m \leq d_m \big) 
	\vee \forall_{n}\ c_n = d_n
$$

\vspace{0.4cm}
\noindent As shown in Lemma \ref{lem:1}  the poset $\Gamma \equiv \langle V, \leq_V \rangle$ is a \textbf{chain}. $\Gamma$ denotes then a new linearly ordered poset comprising  the family  of maximal chains of  $\Pi$ as \textbf{points} of $V$.  The source of this kind of  point of view is inherited from combinatorial interpretation \cite{19}  of $F$-nomial coefficients \cite{10}-\cite{18}, [4,5,6]. Indeed.
Recall \cite{19,18,17} that maximal chains  of layers $\langle \Phi_k \rightarrow \Phi_n \rangle$ are set points of the  set

$$
	\big\{ \mathrm{maximal\ chains\ in\ } \langle \Phi_k \rightarrow \Phi_n \rangle \big\} \equiv
	C_{max}\big( \langle \Phi_k \rightarrow \Phi_n \rangle \big) \equiv
	C^{k,n}_{max}.
$$

\noindent Denoting with $V_{k,n} \subseteq V$ the discrete finite rectangular $F$-hyper-box  or $(k,n)-F$-hyper-box or in everyday parlance just  $(k,n)$-box

$$
	V_{k,n} \equiv [k_F] \times [(k+1)_F] \times ... \times [n_F]
$$

\noindent we identify the two just by agreement  according to  the natural identification:

$$
	C^{k,n}_{max} \equiv V_{k,n}.
	$$

\noindent Here down as in above we shall use [see: Appendix, Section 5] the so called upside down notation: $F_k \equiv  k_F$.

\vspace{0.4cm}
\noindent \textbf{2.2.2} \ $\subseteq$, $B \equiv \langle V, \subseteq \rangle$

\vspace{0.2cm}
\noindent A partial ordered relation $\leq_\Pi = \subseteq$ which sets the pace with the identification  KoDAGs treated as chains of discrete  $F$- hyper-boxes from $N^\infty$ is obvious  and of primary importance - including practical (tilings) and visualization purposes.

\vspace{0.2cm}
The natural partial order $\subseteq$ forced upon discrete boxes $C^{k,n}_{max}$ might be additionally and extra  introduced also in $\Pi =  \langle P, \leq \rangle$  and then it is just  set inclusion. For 

$$
	B \equiv \langle V, \subseteq \rangle = L_{box}(V) \subset L(V) \equiv 2^V
$$

\vspace{0.4cm}
\noindent this would mean correspondingly  the discrete hyper-box into discrete  hyper-box inclusion with the welcomed scenario of an  infinite Boolean lattice $L(V) \equiv 2^V$ or chains of finite Boolean lattices $L(V_{k,n}) \equiv 2^{V_{k,n}}$.

The $\leq_V$  choice is of course different from $\subseteq$ though both are relatives in a sense. The chain poset $\Gamma$ and the sub Lattice $B$  of Boolean Lattice are already  complementarily useful -  both.

More than that. One is not forbidden  -  not at all - to consider simultaneously two partial orders  i.e.  the structure  $\langle V, \leq_\Pi, \subseteq \rangle$ in general  and specifically one of them being linear order $\langle V, \leq_V, \subseteq \rangle$.

\vspace{0.2cm}
\noindent Summing up: \\
\textbf{Mantra:} \textbf{self-quotation:}

\begin{quotation}
\emph{KoDAGs represent the outmost general chains of di-bi-cliques denoting thus the outmost general chains 
of binary relations. Because of  that fact KoDAGs start to become  an outstanding concept of nowadays investigations.}
\end{quotation}

\noindent \textbf{Mantra:} \textbf{Refrazed = coded}:

\begin{quotation}
\emph{KoDAGs in $K \equiv \langle V, \leq_\Pi \rangle$ coding are "just chains" of discrete \textbf{hyper-boxes} thus coding the outmost general
chains  of binary relations. Because of  that fact this KoDAGs in this coding start to become  an outstandingly natural 
concept of nowadays investigations.}
\end{quotation}

\noindent - Hyper-boxes ? \\
\noindent - Yes. For example chains of \emph{di-bi-cliques hyper-boxes}\\
$$V_{k,k+1} = [k_F] \times  [(k+1)_F] \Leftrightarrow C^{k,k+1}_{max}.$$

\noindent Here down as in above we shall use [see: Appendix, Section 5] the so called upside down notation: $F_k \equiv  k_F$.

\vspace{0.6cm}
\noindent \textbf{2.2.3} \ $\leq_\Phi$, $\Phi = \times_{s\geq 0} \Phi_s$ coding ($\Phi$ is infinite Cartesian product of levels $\Phi_s$)

\vspace{0.2cm}
\noindent $\leq_\Pi = \leq_\Phi$. Consider a sequence of posets $\langle\Phi_s, \leq_{\Phi_s} \rangle$ where $\Phi_s$ stays for set of vertices on $s$-th level i.e. $\Phi_s = \big\{ \langle j,s\rangle : 1\leq j \leq s_F \big\}$ and \textbf{linear order} between vertices on this level i.e. 
$$
	\langle j, s\rangle  \leq_{\Phi_s}   \langle k, s\rangle  \Leftrightarrow  j\leq k, \ \ \ \ \  j,k\in [s_F]\equiv \big\{x; \leq x \leq s_F \big\}.
$$

\noindent If so then the  poset $\Lambda = \langle \times_{s\geq 0}\Phi_s, \leq_\Phi \rangle$ is no more cobweb poset  i.e. it is not isomorphic to $\langle P, \leq \rangle$. This is obvious already  from the fact that the sine qua non feature of  $\langle P, \leq \rangle$ is  that $\Phi_s$ are independent sets  and $\langle P_n, \leq \rangle$ is an $n$-level graded poset; $n \in N \cup \{0\}$.

\vspace{0.4cm}
\noindent \textbf{Specifications of $\leq_\Phi$.}
\vspace{0.2cm}

\noindent There are \textbf{three} of the possible "WIKI-common" \textbf{partial orders} on the Cartesian product of \textbf{two} totally ordered sets:

\begin{enumerate}
\item  Lexicographical order: $(a,b) \leq (c,d)$ if and only if $a < c$ or $(a = c \wedge b \leq d)$. This is a \textbf{total order}.

\item \noindent $(a,b) \leq (c,d)$ if and only if $a \leq c$ and $b \leq d$ (the product order). This is a \textbf{partial order}. 

\item \noindent $(a,b) \leq (c,d)$ if and only if $(a < c \wedge b < d)$ or $(a = c \wedge b = d)$ (the reflexive closure of the direct product of the corresponding strict total orders). This results also a \textbf{partial order}.
\end{enumerate}

\noindent As for the $\leq_\Phi$ we have at our disposal unbounded multitude of combining these  three choices at the start   for example for di-bi-cliques and then subsequently imposing various orders on the resulted products of posets  -  freely,  step by step using the induction for example.

\vspace{0.2cm}
\noindent \textbf{Conviction.} 

\noindent Even more ? We are not limited only to this possibility...(?)   [Hypothesis?]

\vspace{0.2cm}
\noindent We leave this for further investigation. Here we content ourselves with an example to be used for visualizations.

\vspace{0.4cm}
\noindent \textbf{Example:} 

\vspace{0.2cm} 

\noindent Consider a poset $\Lambda = \langle \times\Phi_s, \leq_\Phi \rangle$, $s\in N\cup\{0\}$ with Cartesian product of vertices' sets $\Phi = \times_{s\in N\cup\{0\}} = \left\{ (\phi_0, \phi_1, ...) : \phi_n \in \Phi_n \right\}$ and the partial  order $\leq_\Phi$ defined as follows

$$
	\chi_1 \leq_\Phi \chi_2   \Leftrightarrow  \forall_{s\in N\cup\{0\}}\ c_s \leq_{\Phi_s} d_s
$$

\vspace{0.3cm}
\noindent where $\chi_1 = (c_0,c_1,c_2,...), \chi_2 = (d_0,d_1,d_2,...) \in \Phi = \times_{s\in N\cup\{0\}} \Phi_s$

\vspace{0.4cm}
\noindent We obtain a \textbf{lattice} with set points representing  maximal chains in $\Pi$, with partial order relation $\leq_\Phi$. The "level by level" Cartesian product poset $\Lambda$ is a lattice as obviously  for any two  elements from $\Phi$ there exists a supremum (join; $\vee$) and an infimum (meet; $\wedge$).

\noindent \textbf{Now note.} The above choice of $\leq_\Phi$ means  that  both $\chi_1$ and $\chi_2$ are elements of an infinite "vertical"  strip $S_{\chi_1, \chi_2} = \left\{ \chi\in\Phi : \chi_1 \wedge \chi_2 \leq \chi \leq \chi_1 \vee \chi_2 \right\}$ which is also a hyper-box $S_{\chi_1, \chi_2} \subseteq \Phi$. Thus again  discrete hyper-box into discrete  hyper-box inclusion  is inevitably around.   

\vspace{0.4cm}
\noindent \textbf{Illustration}

\begin{figure}[ht]
\begin{center}
	\includegraphics[width=50mm]{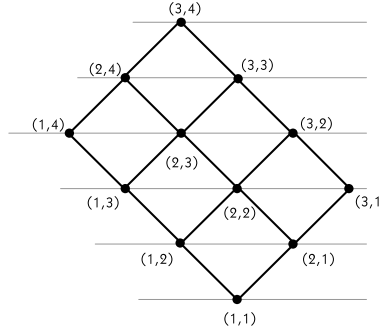}
	\caption{Hasse diagram of lattice $\langle V_{3,4}, \leq_\Phi \rangle$ \label{fig:5}}
\end{center}
\end{figure}

\vspace{0.2cm}
\noindent An illustration of $\Lambda = \langle\Phi, \leq_\Phi \rangle \equiv \langle\times\Phi_s, \leq_\Phi \rangle$ lattice points = maximal chains in  $\Pi$-disposal is supplied by the figure  of the coded di-bi-clique layer $V_{3,4} \subset \Phi$, where $F$ =  Natural numbers (Fig. \ref{fig:5}).

\section{How does it work}

\begin{figure}[ht]
\begin{center}
	\includegraphics[width=70mm]{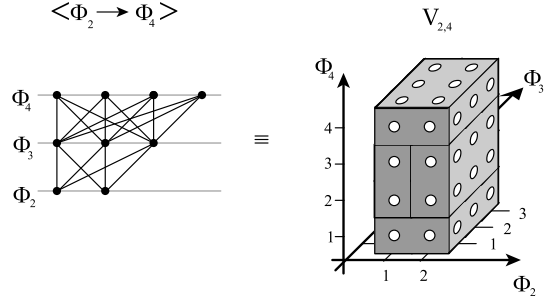}
	\caption{A cobweb layer $\langle\Phi_2 \rightarrow \Phi_4 \rangle$ and equivalent hyper-box $V_{2,4}$ \label{fig:representation}}
\end{center}
\end{figure}

\noindent \textbf{Recall.} Here down as in above we shall use (see: Appendix, Section 5) the so called upside down notation: $F_k \equiv  k_F$.
$V$ is a code of $\Pi$ and $V_{k,n}$ is a code of the layer  $\langle\Phi_k \rightarrow \Phi_n \rangle$. \textbf{Obviously:} any cobweb poset $\Pi$ designated by the sequence $F$ has its own representation $V$ \textbf{independently of the partial order chosen i.e.} $K\equiv\langle V, \leq_\Pi \rangle$. Here comes an illustration (Fig. \ref{fig:representation}).

\noindent As a matter of illustration we show soon using pictures \textbf{how} any cobweb poset $\Pi$ can be represented as a hyper-box $V$ in discrete subspace of $N^\infty$, so that any maximal-chain of $\Pi$ poset becomes  a point of $V$. More than that. 

Employing quite arbitrary sequences $F\equiv\{n_F\}_{n\geq 0}$ infinitely many new representations of natural numbers called an $F$-base or base-$F$ number system representations are introduced. These are used for coding and interpreted as chain coordinatization in KoDAGs as well as systems of  infinite number of  boxes' sequences of $F$-varying containers capacity of subsequent boxes. Needless to say how crucial this base-$F$ number system for KoDAGs - hence - consequently for arbitrary chain s of binary relations. New positional-type number $F$-based systems are umbral base-$F$ number systems in a sense to explained in what follows.

\noindent Let us introduce and recall the basic notions and notation.

\vspace{0.4cm}
\noindent Let $k,n\in N\cup\{0\}$ unless other stated. For $k=n$ we deal then with "empty boxes".

\vspace{0.2cm}
\noindent Denoting with $V_{k,n} \in V$ the discrete finite rectangular $F$-hyper-box or $(k,n)-F$-hyper-box or in everyday parlance just $(k,n)$-box

$$
	V_{k,n} \equiv [k_F]\times [(k+1)_F]\times ... \times[n_F]
$$
\noindent we identify the following two just by agreement according to the $F$-natural identification:
$$
	C^{k,n}_{max} \equiv V_{k,n}
$$

\noindent \textbf{Note.} Any point $\chi \in V_{1,m}$ from finite hyper-box is represented by point from infinite hyper-box $V$ according to the identification
$$
	\chi = \left( c_0,c_1,...,c_m \right) \equiv \left( c_0,c_1,...,c_m,0,0,... \right) \in V 
$$

\begin{defn}[$\lambda$]
We define now $F$-coding map $\lambda:V^{*} \mapsto N\cup\{0\}$ where $V^* = \left\{ (c_0,c_1,...,c_m,0,0,...)\in V\right\}$ simply as 
$$
{
	\lambda\left( (c_0,c_1,...,c_m) \right) = c_0 \cdot 1 + c_1 \cdot k^{\overline{1}}_F + c_2 \cdot k^{\overline{2}}_F + ...
	c_m \cdot k^{\overline{m}}_F
\atop
	\lambda\left(\chi\right) = \sum_{s=0}^{m} c_s \cdot k^{\overline{s}}_F
}
$$
\end{defn}

\noindent where, in Kwa\'sniewski upside down notation (see Section 5)\\
$k^{\overline{s}}_F \equiv k_F(k+1)_F...(k+s-1)_F$, $0^{\overline{0}}_F \equiv 1$ hence $1^{\overline{s}}_F \equiv 1_F(1+1)_F...(1+s-1)_F\equiv s_F!$ hence for $k>0$

$$
	(\lambda)\ \ \ \ \lambda\left((c_0,c_1,...,c_m)\right) = c_0 \cdot 1 +
	\sum_{s=1}^{m}c_s\cdot \frac{(k+s-1)_F!}{(k-1)_F!}
$$
\noindent and consequently for $V_{1,m} \equiv V_m$:

$$
	\lambda\left(\chi\right) = c_0\cdot 1 +
	\sum_{s=1}^{m}c_s\cdot s_F!
$$

\begin{figure}[ht]
\begin{center}
	\includegraphics[width=40mm]{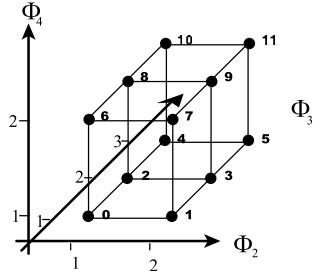}
	\caption{Labeling effect for $V_{2,4}$ \label{fig:labeling}}
\end{center}
\end{figure}

\noindent In not Kwa\'sniewski upside down \textbf{notation}. 

$$
	k^{\overline{s}}_F \equiv k_F(k_F + 1)...(k_F + s-1), \ \ k^{\overline{0}}_F \equiv 1
$$

\noindent hence 
$$
	(\lambda\ not) \ \ \ \  \lambda\left((c_0,c_1,...,c_m)\right) = c_1 \cdot 1 +
	\sum_{s=1}^{m}c_s\cdot \frac{(k_F+s-1)!}{(k_F-1)!}
$$

\noindent and consequently 
$$
	\lambda\left(\chi\right) = c_0\cdot 1 +
	\sum_{s=1}^{m}c_s\cdot s!
$$

\noindent For $F\neq N$, this $\lambda$ map defining formulas in not Kwa\'sniewski upside down notation  above appear  not "$F$-natural" i.e. not $F$-designated structures-consistent. The structures in mind comprise KoDAgs, their incidendece algebras, combinatorial interpretation, tiling problem solutions and so on.

\vspace{0.2cm}
\noindent Let  $F\neq N$.  Then of course

$$
	c_0 \cdot 1 +
	\sum_{s=1}^{m}c_s\cdot \frac{(k+s-1)_F!}{(k-1)_F!}
	\neq
	c_0 \cdot 1 +
	\sum_{s=1}^{m}c_s\cdot \frac{(k_F+s-1)!}{(k_F-1)!}
$$

\noindent It is the Kwa\'sniewski upside down  notation (see Section \ref{sec:appendix})  which appears natural for the important statements to hold and vague "umbral" sequence based positional systems  to work. This is this   very Kwa\'sniewski upside down  notation  to be used here down. 

\vspace{0.2cm}
\noindent For arbitrary admissible $F$ the ($\lambda$) definition formulas in upside down notation  appear "$F$-natural" i.e. $F$-designated structures-consistent. The structures in mind comprise KoDAGs, their incidence algebras, combinatorial interpretation, tiling problem solutions and so on.

\noindent Recall now the notation and the content of the lemma to be used.

\begin{lemma}
$\leq_V$ is total order in $V^*$ hence $\langle V^*, \leq_V \rangle$ is a chain.
\end{lemma}

\begin{lemma}\label{lem:2}
Let $k > 0$. If $\chi_1=(k_F-1, (k+1)_F-1, ..., (k+m-1)_F-1,0) \in V_{k,n}$ and $\chi_2 =(0, 0, ..., 1) \in V_{k,n}$ where $n=k+m \in N$ then 

\begin{equation}\label{eq:lemma}
	\lambda(\chi_1) = \lambda(\chi_2) - 1
\end{equation}
\end{lemma}

\noindent \textit{Proof:} by induction.
\begin{enumerate}
\item For $m=1$ we have $\lambda((k_F-1,0)) = \lambda((0,1)) - 1$, obvious.
\item If (\ref{eq:lemma}) is true for $m>1$ i.e. 
$$
	\sum_{s=0}^{m-1}\left[ (k+s)_F - 1 \right]\cdot k^{\overline{s}}_F = k^{\overline{m}}_F - 1
$$
\noindent then (\ref{eq:lemma}) is true for $(m+1)$ i.e.
$$
	\sum_{s=0}^{m}\left[ (k+s)_F - 1 \right]\cdot k^{\overline{s}}_F = k^{\overline{m+1}}_F - 1.
$$
\end{enumerate}
\noindent Indeed. Just check.\\
$ \sum_{s=0}^{m-1}\left[ (k+s)_F - 1 \right]\cdot k^{\overline{s}}_F + \left[(k+m)_F - 1 \right]\cdot k^{\overline{m}}_F
	= (k+m)_F\cdot k^{\overline{m}}_F - 1$, \\
$ \sum_{s=0}^{m-1}\left[ (k+s)_F - 1 \right]\cdot k^{\overline{s}}_F - k^{\overline{m}}_F + 1 + (k+m)_F\cdot k^{\overline{m}}_F
	= (k+m)_F\cdot k^{\overline{m}}_F$\\
\noindent Then according to the induction assumption $0 + (k+m)_F\cdot k^{\overline{m}}_F = (k+m)_F\cdot k^{\overline{m}}_F$ $\blacksquare$ 

\vspace{0.2cm}
\noindent \textbf{Comment.} This lemma proves a \textbf{$F$-coding numeral system property} via $\lambda$ such as that of the  binary numeral system or base-2 number system as  example shows: $0111_2 = 1000_2 - 1$. 

\vspace{0.2cm}
\noindent Hence the purpose aimed convention $s\in\{0,1,...,(k+m)_F-1\}$ for hyper-box and other definitions in order to exhibit the rule: $11...10 = 00...01$   while moving one step up to the next level of any given KoDAG.

\begin{lemma}\label{lem:3}
For any hyper-box $V_{k,n}$ its $F$-coding, labeling function $\lambda$ is an injection from $V_{k,n}$ to $N\cup\{0\}$ i.e. $\lambda(\chi_1) = \lambda(\chi_2) \Rightarrow \chi_1 = \chi_2$  for any $\chi_1, \chi_2 \in V_{k,n}$.
\end{lemma}

\noindent \textit{Proof} a contrario\\
Let us assume $\lambda(\chi_1) = \lambda(\chi_2) \wedge \chi_1 \neq \chi_2$. Then 
(1) $\sum_{s=0}^{m}c_s\cdot k^{\overline{s}}_F = \sum_{s=0}^{m}d_s\cdot k^{\overline{s}}_F$ $\wedge$ 
(2) $\exists s: c_s \neq d_s$.  
Then $\sum_{s=0}^{m}(c_s - d_s)\cdot k^{\overline{s}}_F = 0$ and from (2) we infer that there exists one or more indices such that 
$(c_s - d_s) \neq 0$. Let $\Omega=\{0,1,...,m\}$ and $I=\{s\in\Omega : (c_s - d_s) \neq 0\}$. Then 
$\sum_{s\in\Omega \setminus I }(c_s - d_s)\cdot k^{\overline{s}}_F = 0$. Therefore we need to consider only
$\sum_{s\in I }(c_s - d_s)\cdot k^{\overline{s}}_F = 0$ case. Let us identify the maximal number $r$ from the set $I$ i.e. 
$r = max\{s\in I\}$. Then (3) 
$\sum_{s\in I\setminus \{r\} }(c_s - d_s)\cdot k^{\overline{s}}_F + (c_s - d_r)\cdot k^{\overline{r}}_F = 0$.
According to Lemma \ref{lem:2} we know that the sum of all $\delta_s\cdot k^{\overline{s}}_F$ where $\delta_s = (c_s - d_s)$, 
$s\in\{0,1,...,r-1\}$ is smaller than $1\cdot k^{\overline{r}}_F$ therefore the sum (3) is not equal to $0$ as seen from 

$$
	\sum_{s\in I\setminus \{r\} }(c_s - d_s)\cdot k^{\overline{s}}_F 
	\leq
	\lambda(\chi_3) < \lambda(\chi_4)
	\leq
	|(c_r - d_r)|\cdot k^{\overline{r}}_F
$$
$$
	\Rightarrow 
	\sum_{s\in I\setminus \{r\} }(c_s - d_s)\cdot k^{\overline{s}}_F + (c_s - d_r)\cdot k^{\overline{r}}_F \neq 0
$$
\noindent where $\chi_3 = (k_F-1, (k+1)_F-1,...,(k+r-1)_F-1,0)$, $\chi_4 = (0,...,0,1)\in V_{k,k+r}$ contrary to the assumption.
$\blacksquare$

\vspace{0.4cm}
\noindent \textbf{Conclusion.} As a conclusion from the above  statements the following is true.

\begin{theoremn}
Chains $\langle V^*,\leq_V \rangle$ and $\langle N, \leq \rangle$ are isomorphic.
\end{theoremn}

\noindent \textit{Proof.}
It is enough to prove that $\lambda$ an order preserving  bijection of $V^*$ and $N$. We already know from Lemma \ref{lem:3} that $\lambda$ is an injection.\\
At first we prove that for any two points $\chi_1, \chi_2 \in V_{k,n}$ where $k,n \in N\cup\{0\}$.
$$
	\chi_1 \leq_V \chi_2 
	\Leftrightarrow 
	\lambda(\chi_1) \leq \lambda(\chi_2)
$$

\noindent Let $\chi_1 = (c_0,c_1,...,c_m)$, $\chi_2 = (d_0,d_1,...,d_m)$. The  case $\chi_1 = \chi_2$ it is obvious, therefore
let us consider $\chi_1 \neq \chi_2$.

\begin{enumerate}
\item $\chi_1 \leq_V \chi_2 \Rightarrow \lambda(\chi_1) \leq \lambda(\chi_2)$

\vspace{0.2cm}
\noindent From definition of $\leq_V$ relation, we have that (1) $\exists_{0\leq j} \leq m (c_j<d_j)$ and (2) $\forall_{t>j} (c_t \leq d_t)$ therefore according to Lemma \ref{lem:2} and (1) we infer 
$\sum_{s=0}^{j}c_s k^{\overline{s}}_F < d_{j} k^{\overline{j}}_F $
and due to (2) 
$$
	\sum_{s=0}^{j}c_s k^{\overline{s}}_F + 
	\sum_{s>j}^{m}c_s k^{\overline{s}}_F < 
	d_{j} k^{\overline{j}}_F +
	\sum_{s>j}^{m}d_s k^{\overline{s}}_F <
	\sum_{s=0}^{j-1}d_s k^{\overline{s}}_F +
	d_{j} k^{\overline{j}}_F +
	\sum_{s>j}^{m}d_s k^{\overline{s}}_F
$$
\noindent Hence $\lambda$ is monomorphism.
\item The $\lambda$ is monomorphism is surjection. This is being proved by supplying an appropriate algorithm  for retrieving $\chi_\alpha$ from its $F$-coding natural number $\alpha$.
\end{enumerate}

\vspace{0.4cm}
\noindent \textbf{Algorithm 1}

\vspace{0.2cm}
\noindent Input: $\alpha\in N\cup\{0\}$, Output: $\chi_\alpha \in V$ such that $\lambda(\chi_\alpha) = \alpha$

\vspace{0.2cm}
\noindent \textbf{Step 1.} 
Find the smallest number $m\in N$ such that $\alpha < k^{\overline{m+1}}_F$, and set $r_m = \alpha$

\vspace{0.2cm}
\noindent \textbf{Step 2.} 
Identify $c_m$ and $r_{m-1}$ from $r_m = c_m\cdot k^{\overline{m}}_F + r_{m-1}$ equation so that $r_{m-1}$ stays for the 
the remainder after appropriate dividing i.e. $r_{m-1} = r_m \mathrm{mod} k^{\overline{m}}_F$ then $r_{m-1} = c_{m-1}\cdot k^{\overline{m-1}}_F + r_{m-2}$ and so on ... $r_1 = c1\cdot k_F + r_0$ until $r_0 = c_0$ or using the while instruction $F$-code algorithm reads:

\vspace{0.3cm}
\noindent \texttt{function lambda($\alpha \in N\cup\{0\}$) : $\chi_\alpha \in V$ \\
begin \\}
{\verb  }\texttt{set $m$ such that $ k^{\overline{m}}_F \leq \alpha < k^{\overline{m+1}}_F$ \\} 
{\verb  }\texttt{$r_m = \alpha$; \\}
{\verb  }\texttt{while ($m > 0$) \\}
{\verb  }\texttt{begin \\}
{\verb         }\texttt{$c_m = r_m$ mod $k^{\overline{m}}_F$; \\}
{\verb         }\texttt{$r_{m-1} = r_m - c_m \cdot k^{\overline{m}}_F$;\\}
{\verb         }\texttt{$m = m - 1$; \\}
{\verb         }\texttt{end \\}
\texttt{$c_0 = r_0$; \\}
\texttt{end }

\vspace{0.4cm}
\noindent This is then just multiple application of Euclidean algorithm and the uniqueness is obvious. 

\vspace{0.2cm}
\noindent \textbf{The result:}

\vspace{0.1cm}
\noindent From Step 1 we have $\alpha = r_m$. From Step 2 $\alpha = c_m\cdot k^{\overline{m}}_F + r_{m-1}$ and so on until $r_0 = c_0$

$$
	\alpha = c_m\cdot k^{\overline{m}}_F + c_{m-1}\cdot k^{\overline{m-1}}_F + ... + c_1\cdot k^{\overline{1}}_F + c_0
	= \lambda(\chi_\alpha) \equiv (c_m c_{m-1} ... c_1 c_0)_F \blacksquare
$$

\noindent \textbf{Name}

\vspace{0.2cm}
\noindent This representation of zero or any natural number $\alpha$ we shall call an \textbf{$F$-base} or \textbf{base-$F$} number system representation of $\alpha$.

\noindent Needless to say how crucial is this \textbf{base-$F$} number system for KoDAGs - hence - consequently for arbitrary chain $s$ of binary relations. 

\vspace{0.4cm}
\noindent \textbf{Examples:}

\vspace{0.2cm}
\noindent Consider now the hyper-box $V_{1,n}$, $n\in N$.

\begin{enumerate}
\item For constant sequence $F\equiv\{p\}_{n \neq 0}$, where $p\in N$ the points of the hyper-box $V(F)$ are $F$-coded with the above use of standard numeral system with the $p$ base or base-$p$ number system i.e. \\
$\lambda\left((c_0,c_1,...,c_m)\right) = (c_m,c_{m-1},...,c_1,c_0)_p = \sum_{s=0}^{m}c_s\cdot p^s$ i.e. \\
$\lambda\left((c_0,c_1,...,c_m)\right) = (c_m c_{m-1} ... c_1 c_0)_p$. \\
\noindent \textbf{Note} the habit-indispensable \textbf{order} of $c_i$ variables \textbf{inversion} in $F$-representation. 
The case $p=2$ is the case of binary system; for example $\chi=(1,0,1,1) \equiv (1101)_2 \equiv (\lambda(\chi))_10 = 11$.
Here then $F$-coding $\lambda$  is a way to convert $F$-base representation to binary numeral base $\lambda((c_0,c_1,...,c_m)) = (c_m c_{m-1}...c_1 c_0)_2$. In general it may be considered as a mean to convert $F$-base representation of a natural number to any other standard base-$d$ representation. For example with $p=10$ choice we arrive at quite ancient decimal system $\lambda\left((c_0,c_1,...,c_m)\right) = (c_m c_{m-1} ... c_1 c_0)_10$ .

\item For non constant sequences like $F$=Fibonacci sequence next weight coefficients are in general different equal to $s_F!$ i.e. \\
$\lambda((c_0,c_1,...,c_m)) = \sum_{s=0}^{m}c_s\cdot k^{\overline{s}}_F$ i.e.
$\lambda((c_0,c_1,...,c_m)) = (c_m c_{m-1}...c_1 c_0)_F$. \\
Recall: the representation  $\lambda((c_0,c_1,...,c_m)) = (c_m c_{m-1}...c_1 c_0)_F$ of zero or any natural number $\alpha$ we call the \textbf{$F$-base} or \textbf{base-$F$} number system representation of $\alpha$.

\vspace{0.2cm}
\noindent Let $F$=Fibonacci sequence. Next weight coefficients are different starting from $s=2$. 
There $\lambda(0,0,0,0,0...) = 0$ ,$\lambda (0,0,1,0,0,0,...)=1\cdot 2_F! = 1$, $\lambda (0,0,0,1,0,0,...) = 1\cdot 3_F! = 2$, 
$\lambda(0,0,1,1,0,0,...)=1\cdot 2_F! + 1\cdot 3_F! = 3$, $\lambda(0,0,0,2,0,0,...) = 2\cdot 3_F! = 4$, 
$\lambda(0,0,1,2,0,0,...)=1\cdot 2_F! + 2\cdot 3_F! = 5$ and so on. For example $32 = \lambda(0,0,0,2,0,1...)= 2\cdot 3_F! + 5_F! = 3_F + 5_F\cdot 4_F \cdot 3_F!$, $24 = \lambda (0,0,0,0,4,0...)= 4\cdot 4F!$, $29 = \lambda (0,0,1,2,4,0...)$, ...

\vspace{0.2cm}
\noindent The way to add one ["ball"] to the store is unique. How it goes? \\
$6 \leftarrow 5+1= \lambda(0,0,1,2,0,0...) + \lambda(0,0,1,0,0,0...) = \lambda(0,0,0,0,1,0...)$ \\
$7 \leftarrow 6+1= \lambda(0,0,0,0,1,0...) + \lambda(0,0,1,0,0,0...) = \lambda(0,0, 1,1,0,0...)$. \\

\vspace{0.2cm}
\noindent The way to add one ["ball"] to the store is unique if one insists to place additional balls one by one.\\
$7 \leftarrow 5+2= \lambda(0,0,1,2,0,0...) + \lambda(0,0,0,1,0,0...) = \lambda(0,0,1,2,0,0...) + \lambda(0,0,1,0,0,0...) + \lambda(0,0,1,0,0,0...)=         \lambda(0,0,0,0,1,0...)+ \lambda(0,0,1,0,0,0...)= \lambda(0,0, 1,1,0,0...)$. This rule is just the $F$-base algorithm of adding:
\begin{verbatim}
  (0,0,1,2,0,0...)
+ (0,0,1,0,0,0...)
------------------
= (0,0,0,0,1,0...)
\end{verbatim}
\end{enumerate}

\noindent The way to add  ball or ball containers or just store to the store  to the store is unique. How it goes?

\noindent We now visualize  this "phenomenon" with boxes, containers and balls, [using \emph{containers} and thus making exact a glimpse's  idea with balls only instead of containers suggested by W. Bajguz - also the participant of our join Gian Carlo Rota Polish Seminar [$http://ii.uwb.edu.pl/akk/sem/sem\_rota.htm$].

\vspace{0.4cm}
\noindent Let $\alpha \in N$. 

\begin{enumerate}
\item \textbf{For constant sequence} $F\equiv\{p\}_{n\geq 0}$ imagine infinite number of  boxes' sequence of the $(p-1)$ containers capacity - each box, as Fig. \ref{fig:box1} shows. Almost all boxes are empty. So containers contain containers. And ultimately - containers in containers which are in containers and so on - store balls.

\begin{figure}[ht]
\begin{center}
	\includegraphics[width=80mm]{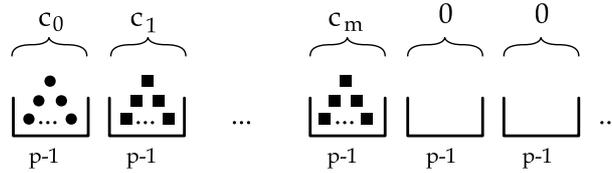}
	\caption{ Containers' sequence interpretation  of the p-base numeral system \label{fig:box1}}
\end{center}
\end{figure}

\vspace{0.2cm}
\noindent \textbf{The rule} of the base-$F$  (i.e. here base-$p$) number system is the following: 
\emph{"adding $p$-th container to $s$-th box causes the surplus and this whole surfeit is passed over as the one container to the next $(s+1)$-th box"} leaving the $s$-th box void of containers. 

\vspace{0.2cm}
\noindent Delivering $\alpha$ balls one by one into $0$-th box and applying \textbf{The rule} one arrives at the final distribution of containers of containers etc. with balls. Distribution is then just Russian \emph{"Babushkas in Babushka"} sequence. Naturally application of this incrassating \emph{"Babushka in Babushka"} rule i.e. The rule ends up with no more than $(p-1)$ containers per box distribution (Figure \ref{fig:box1}).

\noindent In zero-th box there are $c_0$ balls $\bullet ...\bullet$. In $s$-th box $s>0$, there are $c_s$ containers $\blacksquare...\blacksquare$.

\item \textbf{For non constant sequence} $F\equiv\{p\}_{n\geq 0}$ \\
\textbf{The F-rule} of any  base-$F$  number system in which the numerical base varies from position to position is the following: \\
\emph{"adding  one $s_F$-th ball to the s-th box causes the surplus and this whole surfeit is passed over as \textbf{the one} container to the next $(s+1)$-th box leaving the $s$-th box void of containers."}

\vspace{0.2cm}
\noindent Another words: the $s_F$ containers overflow in the $s$-th box is passed over as an additional \textbf{one container} into the $(s+1)$ box $s$-th box void of containers.\\

\noindent Naturally application of this incrassating different in size \emph{"Babushka in Babushka"} rule i.e. \textbf{The F-rule} ends up with no more than $(s_F-1)$ containers per $s$-th box distribution (Figure \ref{fig:box2}).

\begin{figure}[ht]
\begin{center}
	\includegraphics[width=80mm]{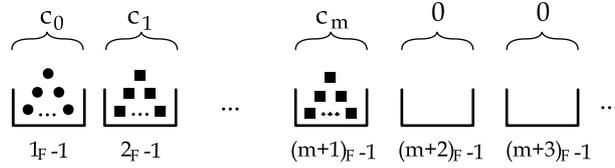}
	\caption{Sequence of containers' interpretation  of the F-base numeral system \label{fig:box2}}
\end{center}
\end{figure}

\noindent In zero-th box there are $c_0$ balls $\bullet ...\bullet$. In $s$-th box $s>0$, there are $c_s$ containers $\blacksquare...\blacksquare$.
\end{enumerate}

\vspace{0.2cm}
\noindent Recall: $c=(c_0,c_1,...,c_m,0,0,...)$. Observe: the number of containers in the $s$-th box equals to the $s$-th coordinate of the point $\chi\in V(F)$

\vspace{0.4cm}
\noindent \textbf{Summing up:}

$$
	\alpha = c_m\cdot k^{\overline{m}}_F + c_{m-1}\cdot k^{\overline{m-1}}_F + ... + c_1\cdot k^{\overline{1}}_F + c_0 = \lambda(\chi_\alpha)
$$
$$
	\lambda\left((c_0,c_1,...,c_m)\right) \equiv (c_m c_{m-1} ... c_1 c_0)_F
$$

\noindent This representation $\lambda\left((c_0,c_1,...,c_m)\right) \equiv (c_m c_{m-1} ... c_1 c_0)_F$  of zero or any natural number $\alpha$ we shall call an \textbf{$F$-base} or \textbf{base-$F$} number system representation of $\alpha$.

\noindent Needless to say how crucial is this \textbf{base-$F$} number system for KoDAGs - hence - consequently for arbitrary chain s of binary relations. 

\noindent New positional-type number $F$-based systems are umbral \textbf{base-$F$} number systems in a following sense explained symbolically in a pictogram sequence way

$$
	\sum_{s=0}^{m} c_s\cdot a^s \rightarrow \mathrm{[umbral\ view]} \rightarrow 
	\sum_{s=0}^{m} c_s\cdot a_s \rightarrow 
$$
$$
	\rightarrow \mathrm{[example\ and\ upside\ down\ notation]} \rightarrow 
	\sum_{s=0}^{m} c_s\cdot k^{\overline{s}}_F \rightarrow 
$$
$$
	\rightarrow \mathrm{[example's\ case\ } k=1, \mathrm{ akk\ upside-down\ notation ]} \rightarrow
	\sum_{s=0}^{m} c_s\cdot s_F! 
$$

\vspace{0.4cm}
\noindent The base-$F$ number system is to be afterwards and soon compared with using the Fibonacci numbers to represent whole numbers\\
$[http://ii.uwb.edu.pl/akk/suprasl/akk.htm]$.

\noindent We mean by this the Fibonacci base system. Recall: \textbf{Zeckendorf } proved in 1972  \cite{31} that:
each representation of a number $n$ as a sum of distinct Fibonacci numbers, is unique but where no two consecutive Fibonacci numbers are used (and there is only one column headed "1").  Another words, a number $\alpha$ might be written as a sum of nonconsecutive Fibonacci numbers $\alpha = \sum_{s=0}^{m}c_s F_s$ where $c_s$ are $0$ or $1$ and $c_s\cdot c_{s+1} = 0$

\section{Application. Cobweb Hyper-Box Tillig Phenomenon.}

\noindent For Kwa\'sniewski tiling problem solution see \cite{2,3}.

\vspace{0.2cm} 
\noindent Here we refraze the tile notion in Dziemia\'nczuk \textbf{geometric-coded} setting in order to accomplish illustrative pictures' delivery. 

\vspace{0.4cm}
\noindent \textbf{Geometric code description of the \emph{cobweb tile}.}

\vspace{0.4cm}
\noindent The set $\tau_m$ of $\lambda$ points $\pi \in V_{k,n}$ and a permutation $\sigma$ of the set $\{1_F,2_F,...,m_F\}$ such that

$$
	\tau_m \equiv \big\{ \pi = (c_1,...,c_m) : c_s \in [(\sigma\cdot s)_F], s=1,2,...,m \big\}
$$

\vspace{0.3cm}
\noindent where $\lambda = m_F!, m=n-k+1$ is called a \emph{cobweb tile}.

\begin{figure}[ht]
\begin{center}
	\includegraphics[width=90mm]{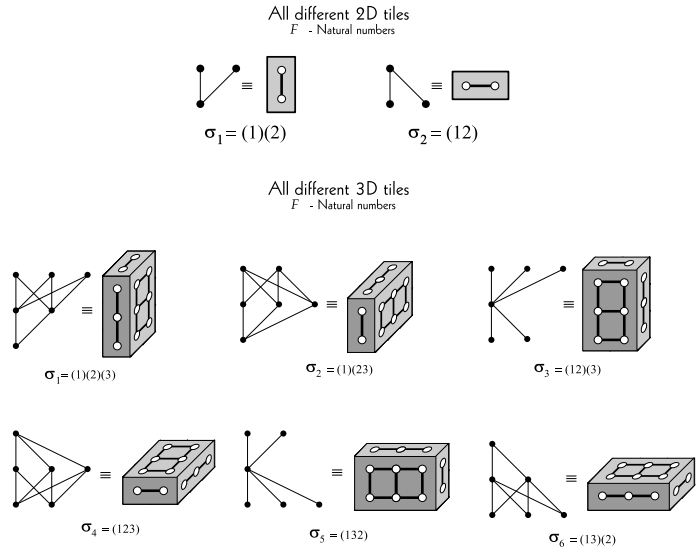}
	\caption{Picture of all Natural numbers' tiles $\tau_2, \tau_3$. \label{fig:tilesN}}
\end{center}
\end{figure}

\begin{figure}[ht!]
\begin{center}
	\includegraphics[width=90mm]{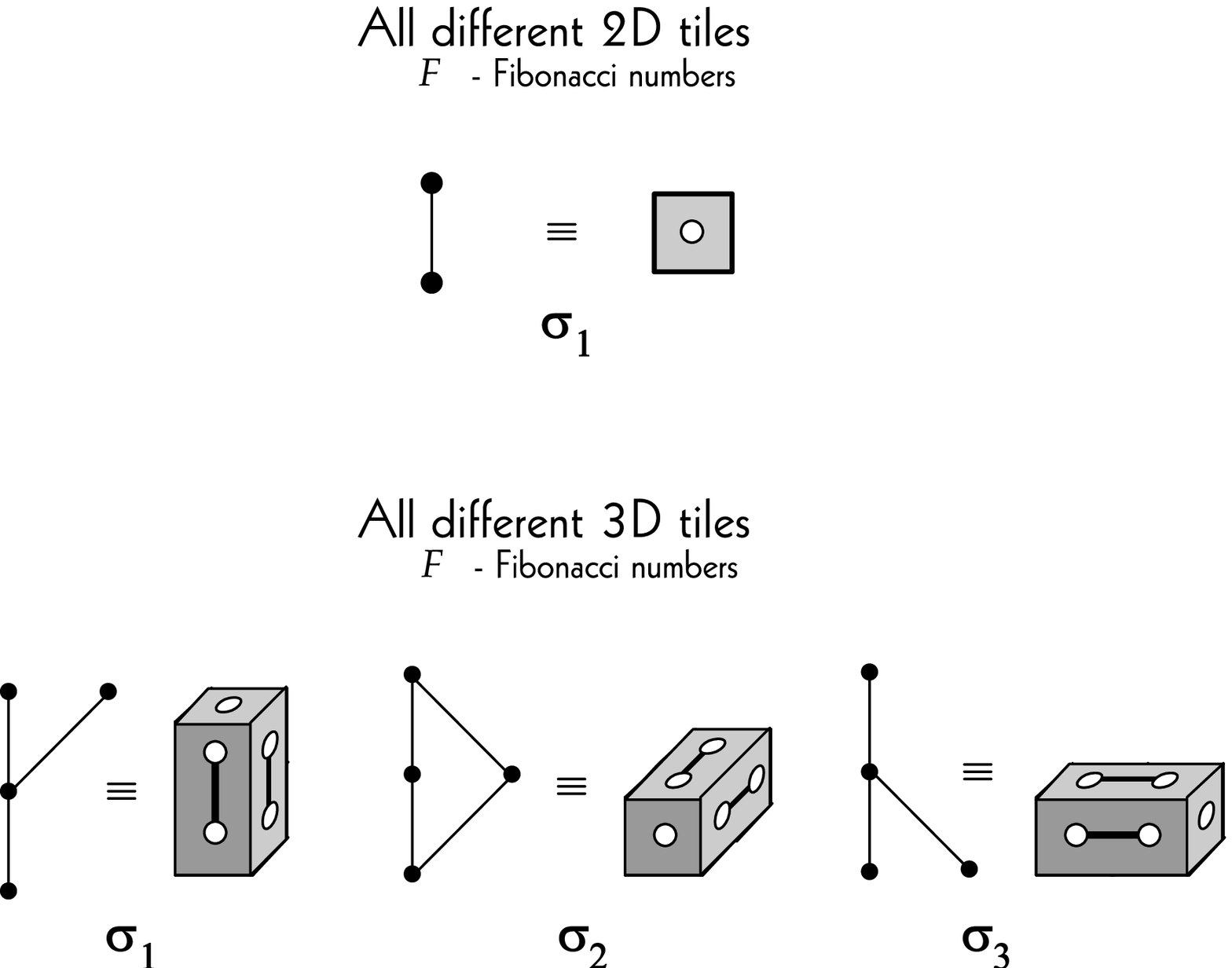}
	\caption{Picture of all Fibonacci numbers' tiles $\tau_2, \tau_3$. \label{fig:tilesF}}
\end{center}
\end{figure}

\vspace{0.4cm}
\noindent Here come some pictures of cobweb hyper-box tiles below (Fig. \ref{fig:tilesN}, \ref{fig:tilesF}).

\begin{figure}[ht!]
\begin{center}
	\includegraphics[width=90mm]{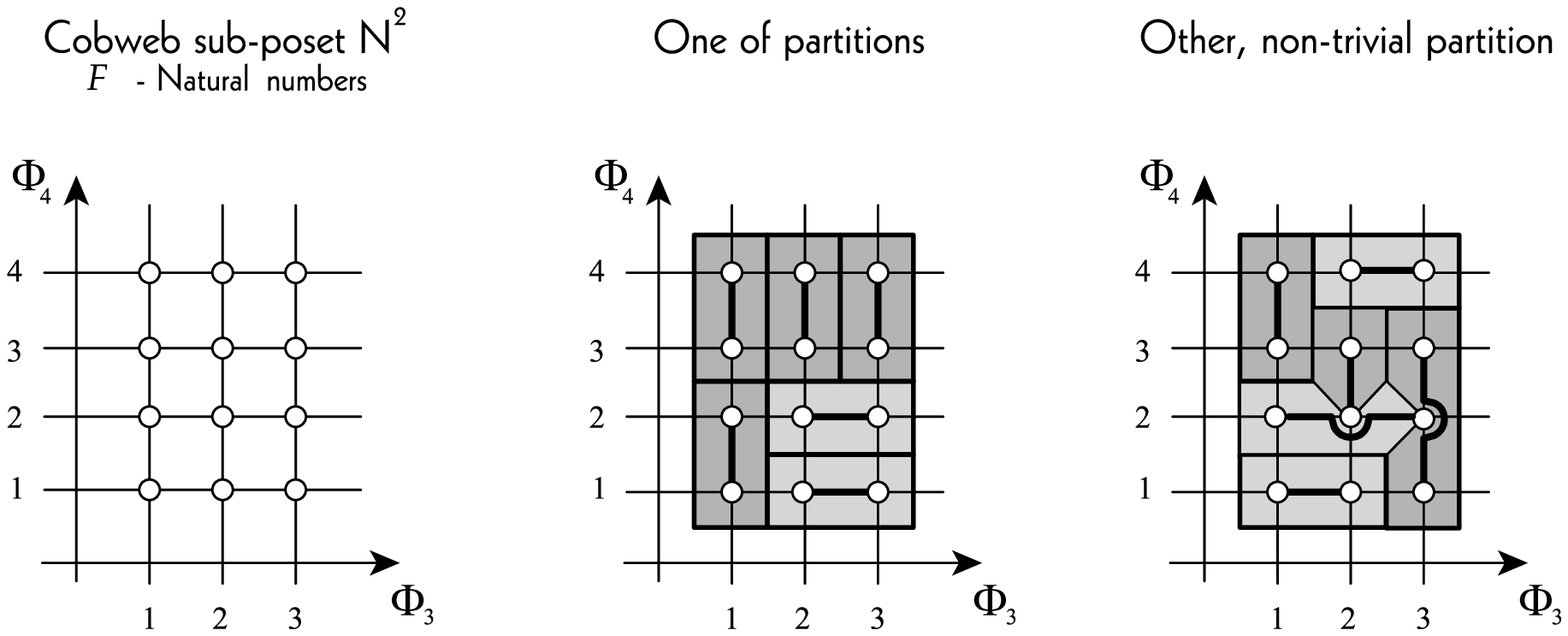}
	\caption{Sample tilings of Natural numbers' $V_{3,4}$. \label{fig:part2N}}
\end{center}
\end{figure}

\begin{figure}[ht!]
\begin{center}
	\includegraphics[width=70mm]{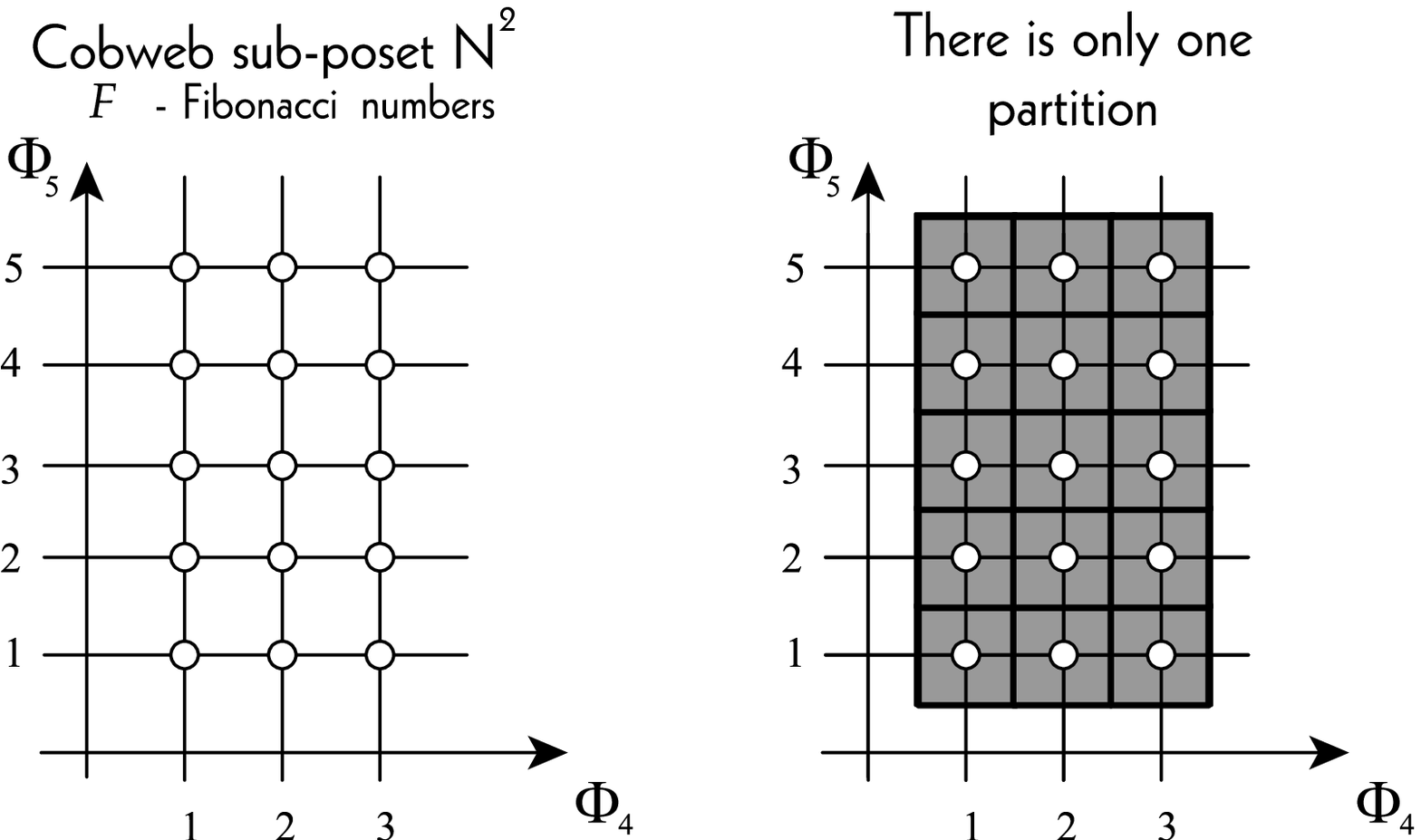}
	\caption{Sample tiling of Fibonacci numbers' $V_{4,5}$. \label{fig:part2F}}
\end{center}
\end{figure}

\begin{figure}[ht!]
\begin{center}
	\includegraphics[width=90mm]{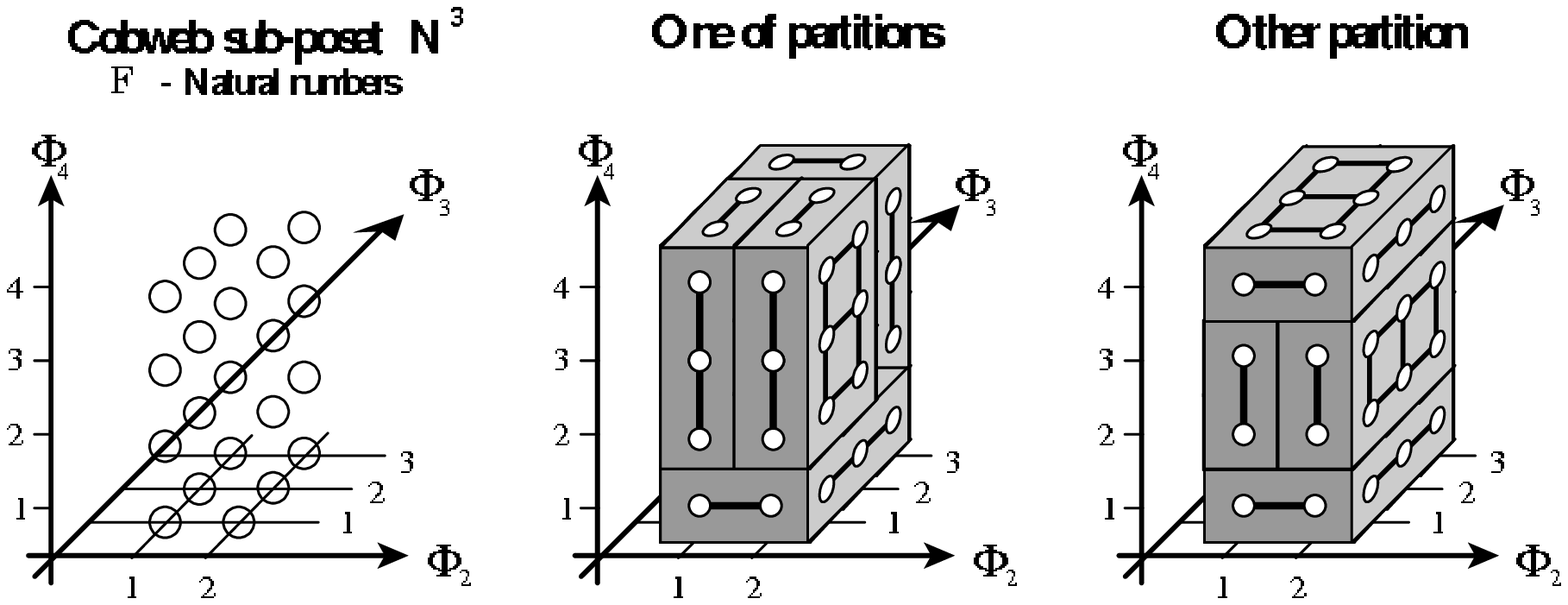}
	\caption{Sample tilings of Natural numbers' $V_{2,4}$. \label{fig:part3N}}
\end{center}
\end{figure}

\begin{figure}[ht!]
\begin{center}
	\includegraphics[width=90mm]{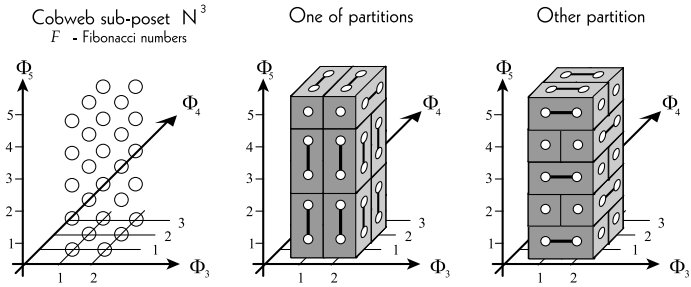}
	\caption{Sample tilings of Fibonacci numbers' $V_{3,5}$. \label{fig:part3F}}
\end{center}
\end{figure}

\vspace{0.4cm}
\noindent \textbf{Definition of \emph{cobweb tiling} in the geometric code.}

\vspace{0.4cm}
\noindent A set $T_{k,n}$ of tiles $\tau_s$ from finite cobweb hyper-box $V_{k,n}$ such that

$$
	T_{k,n} = \big\{ \tau_s \in V_{k,n} : \tau_i \cap \tau_j = \emptyset, i\neq j \wedge \bigcup_{i=1}^{m_F!}\tau_i = V_{k,n} \big\}
$$

\noindent is called \emph{cobweb tiling}.

\vspace{0.4cm}
\noindent \textbf{Notation.} Let $\mathcal{T}(V_{k,n})$ denotes the family of all tilings $T_{k,n}$ of finite cobweb hyper box $V_{k,n}$.

\vspace{0.4cm}
\noindent \textbf{Information.} The compact formula for the number $\left\{{n \atop k}\right\}_F$ of all tilings of the layer $\langle\Phi_k\rightarrow\Phi_n\rangle$ i.e. $\left\{{n \atop k}\right\}_F = |\mathcal{T}(V_{k,n})|$ is still not known \cite{19,17}. See more in \cite{28,2,3}.

\vspace{0.4cm}
\noindent To this end by now we supply  some pictures of the cobweb hyper-box tilings  with (Fig. \ref{fig:part2N} - \ref{fig:part3F})

\vspace{0.4cm}
\noindent \textbf{To be continued.}

\vspace{0.4cm}
\noindent \textbf{Summing up.}

\vspace{0.2cm}
\noindent Cobweb hyper-box representation of cobweb poset in infinite discrete space as a rectangular hyper box invented as a result of cobweb tilings' visualization  [\emph{http://www.dejaview.cad.pl/cobwebposets.html}]  
already shows up efficient and promising for future investigations. Further readings:  \cite{19,18,17,1,2,3,28}.


\section{Appendix - Kwa\'sniewski upside down  notation }\label{sec:appendix} 

\noindent  The Kwa\'sniewski  upside - down notation $n_F \equiv F_n$  has been used for mnemonic reasons - as in the case of Gaussian numbers in finite geometries  and the so called \emph{"quantum groups"}; (see [10,11], [18-20] and references therein). It has been used then consequently in all relevant papers by Ewa Krot, by Ewa Krot-Sieniawska and Dziemia\'nczuk and Bajguz.

\vspace{0.2cm}

\noindent Precise formulation of the Upside Down Principle may be found in  recent  (20 Feb 2009) The Internet Gian Carlo Rota Polish Semenar affiliated article  [32]. 
There one also finds the recent formulation of Kwa\'sniewski combinatorial interpretation of the $F$-nomial coefficients (consult Appendix 2. for more on that).

\vspace{0.2cm}
Given any sequence $\{F_n\}_{n\geq 0}$ of nonzero reals
($F_0 = 0$ being sometimes acceptable as $0! = F_0! = 1.$)
one defines  its corresponding  binomial-like $F-nomial$  coefficients as in Ward`s Calculus of sequences [27] as follows.
\begin{defn}.$(n_{F}\equiv F_{n}\neq 0,\quad n > 0)$
$$
\left( \begin{array}{c} n\\k\end{array}
\right)_{F}=\frac{F_{n}!}{F_{k}!F_{n-k}!}\equiv
\frac{n_{F}^{\underline{k}}}{k_{F}!},\quad \quad n_{F}!\equiv n_{F}(n-1)_{F}(n-2)_{F}(n-3)_{F}\ldots 2_{F}1_{F};$$
$$ 0_{F}!=1;\quad n_{F}^{\underline{k}}=n_{F}(n-1)_{F}\ldots (n-k+1)_{F}. $$
\end{defn}

\noindent  Kwa\'sniewski  hade made above  an analogy driven identifications in the spirit of  Ward`s Calculus of sequences [27]. 
Identification $n_{F}\equiv F_{n}$ is  the notation used in extended Fibonomial Calculus case [10,11-15,4,5,6] being also there inspiring as $n_F$ mimics $n_q$ established notation for Gaussian integers exploited in much elaborated family of various applications including quantum physics (see [10,11] and references therein).

\vspace{0.2cm}
\noindent \textbf{Now compare:  [10,11]}

$$
	{n \choose k} = \frac{n!}{k!(n-k)!} \rightarrow
	{n \choose k}_q = \frac{n_q!}{k_q!(n-k)_q!} \rightarrow
	{n \choose k}_\psi = \frac{n_\psi!}{k_\psi!(n-k)_\psi!}
$$

\noindent 
${n \choose k} \rightarrow {n \choose k}_q \rightarrow {n \choose k}_\psi$,... incidence-like coefficients, connection constants 
\textbf{see: [25,26]} and references therein; \textbf{here and there throughout up-side down notation is mnemonic source of  associations}

$$
	\bigg\{{n \atop k }\bigg\} \rightarrow \bigg\{{n \atop k }\bigg\}_q \rightarrow \bigg\{{n \atop k }\bigg\}_\psi, ... 
	\ \ \ \
	\bigg[{n \atop k }\bigg] \rightarrow \bigg[{n \atop k }\bigg]_q \rightarrow \bigg[{n \atop k }\bigg]_\psi, ... 
$$

\vspace{0.2cm}
\noindent \textbf{$\psi$ denotes an extension of}
$$
	\bigg\{\frac{1}{n!}\bigg\}_{n\geq 0}
$$

\noindent sequence to  quite arbitrary one (\emph{"admissible"}) and the specific choices are for example: Fibonomialy-extended ($F_n$, $n\geq 0$ -Fibonacci sequence) or Gauss $q$-extended 

$$
	\left\{ \psi_n \right\}_{n\geq 0} = \bigg\{\frac{1}{F_n!}\bigg\}_{n\geq 0}, \ \ \ \ 
	\left\{ \psi_n \right\}_{n\geq 0} = \bigg\{\frac{1}{n_q!}\bigg\}_{n\geq 0}, \ \ \ \ 
$$

\noindent admissible  sequences of extended umbral operator calculus - see more below. With such an extension we may $\psi$-mnemonic repeat with exactly the same simplicity and beauty much of what was done by Rota years ago. Thus via practisizing  we get used to write down these extensions in mnemonic upside down notation [1-25]:

$$
	n_\psi \equiv \psi_n, \ \ x_\psi \equiv \psi(x) \equiv \psi_x, \ \ n_\psi! = n_\psi(n-1)_\psi!, \ \ 0_\psi = 1
$$
$$
	x^{\underline{k}}_\psi = x_\psi(x-1)_\psi...(x-k+1)_\psi \equiv \psi(x)\psi(x-1)...\psi(x-k+1)
$$

\vspace{0.2cm}
\noindent You may consult for further development and use of this notation [10,11,22] and references therein.  As for references - the papers of main references are: [10,11,22]. 

\vspace{0.2cm}
 
\noindent The Kwa\'sniewski upside-down notation imposes associations with more general schemes notions
as  Whitney numbers  aside of incidence coefficients   [18]                  

\vspace{0.2cm}
\noindent \textbf{Summarizing}. 
While identifying general properties of such $\psi$-extensions in all their connotations mentioned above the merit consists indeed in notation i.e. here - in writing objects of these extensions in mnemonic convenient \textbf{upside down notation} .   
\begin{equation}\label{eq31}
\frac {\psi_{(n-1)}}{\psi_n}\equiv n_\psi,
n_\psi!=n_\psi(n-1)_\psi!, n>0 ,   x_{\psi}\equiv \frac
{\psi{(x-1)}}{\psi(x)} ,
\end{equation}
\begin{equation}\label{eq32}
x_{\psi}^{\underline{k}}=x_{\psi}(x-1)_\psi(x-2)_{\psi}...(x-k+1)_{\psi}
\end{equation}
\begin{equation}\label{eq33}
x_{\psi}(x-1)_{\psi}...(x-k+1)_{\psi}=
\frac{\psi(x-1)\psi(x-2)...\psi(x-k)} {\psi(x)
\psi(x-1)...\psi(x-k +1)} .
\end{equation}
If one writes the above in the form $x_{\psi} \equiv \frac
{\psi{(x-1)}}{\psi(x)}\equiv \Phi(x)\equiv\Phi_x\equiv x_{\Phi}$ ,
\noindent one sees that the name upside down notation is legitimate.

\vspace{0.5cm}

\noindent \textbf{ The KoDAG Enterprise Information.}  For most recent (20 Feb 2009) developments on upside down notation efficiency see  [32]   and 

\vspace{0.2cm}

\textit{This is The Definition  which is the right answer to the leitmotiv questions}

\vspace{0.1cm}

\noindent i.e. \textcolor{green}{\textbf{The Internet Gian-Carlo Polish Seminar article}},\\
No \textbf{1}, \textbf{Subject 2}, \textbf{2009-02-18}\\
\noindent \emph{http://ii.uwb.edu.pl/akk/sem/sem\_rota.htm}\\

\section{Appendix - Cobweb posets and  KoDAGs' ponderables  of Kwa\'sniewski relevant recent productions. Excerpts from [32].}\label{sec:appendix}

\textbf{A.2.1}

\begin{defn}
\noindent  Let  $n\in N \cup \left\{0\right\}\cup \left\{\infty\right\}$. Let   $r,s \in N \cup \left\{0\right\}$.  Let  $\Pi_n$ be the graded partial ordered set (poset) i.e. $\Pi_n = (\Phi_n,\leq)= ( \bigcup_{k=0}^n \Phi_k ,\leq)$ and $\left\langle \Phi_k \right\rangle_{k=0}^n$ constitutes ordered partition of $\Pi_n$. A graded poset   $\Pi_n$  with finite set of minimal 
elements is called \textbf{cobweb poset} \textsl{iff}  
$$\forall x,y \in \Phi \  i.e. \  x \in \Phi_r \ and \  y \in \Phi_s \   r \neq s\ \Rightarrow \   x\leq y   \ or \ y\leq x  , $$ 
 $\Pi_\infty \equiv \Pi. $
\end{defn}

\vspace{0.1cm}

\noindent \textbf{Note}. By definition of $\Pi$ being graded its  levels    $\Phi_r \in \left\{\Phi_k\right\}_k^\infty$ are independence sets  and of course partial order  $\leq $ up there in Definition 6.1. might be replaced by $<$.

\vspace{0.2cm}

\noindent The Definition ?.  is the reason for calling Hasse digraph $D = \left\langle \Phi, \leq \cdot \right\rangle $ of the poset $(\Phi,\leq))$ a \textbf{\textcolor{red}{Ko}}DAG as in  Professor   
\textbf{\textcolor{red}{K}}azimierz   \textbf{\textcolor{red}{K}}uratowski native language one word \textbf{\textcolor{red}{Ko}mplet} means \textbf{complete ensemble} - see more in  [23]
and for the history of this name see:  The Internet Gian-Carlo Polish Seminar \textbf{Subject 1.  oDAGs and KoDAGs in Company} (Dec. 2008).

\begin{defn}
\noindent Let  $F = \left\langle k_F \right\rangle_{k=0}^n$ be an arbitrary natural numbers valued sequence, where $n\in N \cup \left\{0\right\}\cup \left\{\infty\right\}$. We say that the cobweb poset $\Pi = (\Phi,\leq)$ is \textcolor{red}{\textbf{denominated}} (encoded=labelled) by  $F$  iff   $\left|\Phi_k\right| = k_F$ for $k = 0,1,..., n.$
\end{defn}

\noindent \textbf{A.2.2}

\begin{observen}
 Let  $n = k+m $. The number of subposets equipotent to subposet
$P_{m}$ rooted at any \textbf{fixed}  point  at the level labeled by
$F_{k}$ and ending at the $n$-th level  labeled by  $F_{n}$ is
equal to
$$ \fnomial{n}{m}=\fnomial{n}{k}=\frac{n^{\underline{k}}_{F}}{k_{F}!}.$$
\end{observen}

\vspace{0.2cm}

\noindent Equivalently - now in a bit more mature \textbf{2009} year the  answer  is given simultaneously viewing layers as biunivoquely representing maximal chains sets. Let us make it formal.

\noindent Such recent equivalent formulation of this combinatorial interpretation is to be found in [33] from where we quote it here down  (see also [32]).

\noindent Let $\left\{ F_n \right\}_{n\geq 0}$ be a natural numbers valued sequence with $F_0 = 1$ (or  $F_0! \equiv 0!$ being exceptional as in case of Fibonacci numbers). Any such sequence uniquely designates both $F$-nomial coefficients of an $F$-extended umbral calculus as well as $F$-cobweb poset introduced  by this author (see :the source [19] from 2005 and earlier references therein). If these $F$-nomial coefficients are natural numbers or zero then we call the sequence $F$ - the $F$-\textbf{cobweb admissible sequence}.

\vspace{0.2cm}

\begin{defn}
Let any $F$-cobweb admissible sequence be given then $F$-nomial coefficients are defined as follows
$$
	\fnomial{n}{k} = \frac{n_F!}{k_F!(n-k)_F!} 
	= \frac{n_F\cdot(n-1)_F\cdot ...\cdot(n-k+1)_F}{1_F\cdot 2_F\cdot ... \cdot k_F}
	= \frac{n^{\underline{k}}_F}{k_F!}
$$
\noindent while $n,k\in N $ and $0_F! = n^{\underline{0}}_F = 1$.
\end{defn}

\vspace{0.2cm}

\begin{defn}

$C_{max}(P_n) \equiv  \left\{c=<x_0,x_1,...,x_n>, \: x_s \in \Phi_s, \:s=0,...,n \right\} $ i.e. $C_{max}(P_n)$ is the set of all maximal chains of $P_n$
\end{defn}

\vspace{0.2cm}

\begin{defn}
Let  $$C_{max}\langle\Phi_k \to \Phi_n \rangle \equiv \left\{c=<x_k,x_{k+1},...,x_n>, \: x_s \in \Phi_s, \:s=k,...,n \right\}.$$
Then the $C\langle\Phi_k \to \Phi_n \rangle $ set of  Hasse sub-diagram corresponding maximal chains defines biunivoquely 
the layer $\langle\Phi_k \to \Phi_n \rangle = \bigcup_{s=k}^n\Phi_s$  as the set of maximal chains' nodes and vice versa -
for  these \textbf{graded} DAGs (KoDAGs included).
\end{defn}

\vspace{0.2cm}

\noindent The \textbf{ equivalent} to those from [17,19]  formulation of combinatorial interpretation of cobweb posets via their cover relation digraphs (Hasse diagrams) is the following.

\vspace{0.2cm}

\noindent \textbf{Theorem} [33]\\
\noindent(Kwa\'sniewski) \textit{For $F$-cobweb admissible sequences $F$-nomial coefficient $\fnomial{n}{k}$ is the cardinality of the family of \emph{equipotent} to  $C_{max}(P_m)$ mutually disjoint maximal chains sets, all together \textbf{partitioning } the set of maximal chains  $C_{max}\langle\Phi_{k+1} \to \Phi_n \rangle$  of the layer   $\langle\Phi_{k+1} \to \Phi_n \rangle$, where $m=n-k$.}

\vspace{0.2cm}

\vspace{0.2cm}
\noindent For February 2009 readings on further progress in \textbf{combinatorial interpretation and application} of the  partial order sets named  cobweb posets and their's corresponding encoding Hasse diagrams KoDAGs see  [32-37] and references therein. For active presentation of cobweb posets see [38].


\vspace{0.3cm} 

\noindent \textbf{Acknowledgments} The authors appraise dr W. Bajguz  and dr Ewa Krot-Sieniawska for the interest in this investigation

\end{document}